# EXTENSION OF A SUMMATION DUE TO RAMANUJAN

ARJUN K. RATHIE

Department of Mathematics, School of Mathematical and Physical Sciences,

Central University of Kerala, Riverside Transit Campus, Padennakad,

P.O. Nileshwar, KASARAGOD, 671328, KERALA, INDIA

e-mail address: akrathie@cukerala.edu.in

**Abstract :** In this short research note, we aim to establish an interesting extension of a summation due to Ramanujan. The result is derived with the help of extension of Gauss summation theorem available in the literature.

## 1. INTRODUCTION

The generalized hypergeometric function $_pF_q$ with p numeratorial and q denominatorial parameters is defined by [7],

$$_pF_q \left[ \begin{matrix} a_1, \ldots, a_p \\ b_1, \ldots, b_q \end{matrix} ; z \right] = \sum_{n=0}^{\infty} \frac{(a_1)_n \ldots (a_p)_n}{(b_1)_n \ldots (b_q)_n} \frac{z^n}{n!} \qquad (1.1)$$

where $(a)_n$ denotes the Pochhammer symbol (or the shifted or raising factorial, since $(1)_n = n!$) defined for any complex number $a$ by

$$(a)_n = \begin{cases} 1, & n = 0 \\ a(a+1)\ldots(a+n-1), & n \in N \end{cases} \qquad (1.2)$$

Using the fundamental function relation, $\Gamma(a+1) = a\,\Gamma(a)$, $(a)_n$ can be written as

$$(a)_n = \frac{\Gamma(a+n)}{\Gamma(a)}, \qquad (n \in N \cup \{0\}) \qquad (1.3)$$

where $\Gamma$ is the well known Gamma function.

For convergence of ${}_pF_q[z]$ and its other elementary properties, we refer [7].

It is interesting to mention here two special cases of (1.1) for (i) $p = 2, q = 1$ and for (ii) $p = 1, q = 1$ respectively defined by

$${}_2F_1\left[\begin{matrix}a, b\\c\end{matrix}; z\right] = \sum_{n=0}^{\infty} \frac{(a)_n (b)_n}{(c)_n} \frac{z^n}{n!} \tag{1.4}$$

and

$${}_1F_1\left[\begin{matrix}a\\c\end{matrix}; z\right] = \sum_{n=0}^{\infty} \frac{(a)_n}{(c)_n} \frac{z^n}{n!} \tag{1.5}$$

The function ${}_2F_1$ is popularly known as the Gauss's hypergeometric function and the function ${}_1F_1$ is known as Kummer's function or the confluent hypergeometric function. The functions ${}_2F_1$ and ${}_1F_1$ form the core of special functions and include as their special cases almost all of the commonly used elementary functions.

It is well known that in the theory of hypergeometric and generalized hypergeometric series, classical summation theorem such as those of Gauss, Gauss second, Kummer and Bailey for the series ${}_2F_1$; Watson, Dixon, Whipple and Saalschutz for the series ${}_3F_2$ and others play an important role. Here we would like to mention two of them. These are:

**Gauss's summation theorem [7]:**

$${}_2F_1\left[\begin{matrix}a, b\\c\end{matrix}; 1\right] = \frac{\Gamma(c)\,\Gamma(c-a-b)}{\Gamma(c-a)\,\Gamma(c-b)} \tag{1.6}$$

provided $\Re(c - a - b) > 0$.

and

**Dixon's summation theorem [7] :**

$$_3F_2\left[\begin{array}{ccc} a, & b, & c \\ 1+a-b, & & 1+a-c \end{array}; 1\right]$$

$$= \frac{\Gamma\left(1+\frac{a}{2}\right)\Gamma(1+a-b)\Gamma(1+a-c)\Gamma\left(1+\frac{a}{2}-b-c\right)}{\Gamma(1+a)\Gamma\left(1+\frac{a}{2}-b\right)\Gamma\left(1+\frac{a}{2}-c\right)\Gamma(1+a-b-c)} \qquad (1.7)$$

provided $\Re(a-2b-2c) > -2$

Recently good deal of progress has been done in generalizing and extending [2,4] the above mentioned classical summation theorems [2,4] for the series $_2F_1$, $_3F_2$ and others. Here we will mention the extensions of (1.6) and (1.7) given in [2,3]

$$_3F_2\left[\begin{array}{ccc} a, & b, & d+1 \\ & c+1, & d \end{array}; 1\right]$$

$$= \frac{\Gamma(c+1)\,\Gamma(c-a-b)}{\Gamma(c-a+1)\,\Gamma(c-b+1)}\left[(c-a-b)+\frac{ab}{d}\right] \qquad (1.8)$$

provided $\Re(c-a-b) > 0$ and $d \neq 0, -1, -2, \cdots$

and

$$_4F_3\left[\begin{array}{cccc} a, & b, & c, & d+1 \\ 2+a-b, & 1+a-c, & d \end{array}; 1\right]$$

$$= \frac{\alpha}{(b-1)} \frac{2^{-a}\Gamma(1+a-c)\,\Gamma(2+a-b)\,\Gamma\left(\frac{3}{2}+\frac{1}{2}a-b-c\right)\Gamma\left(\frac{1}{2}\right)}{\Gamma\left(\frac{1}{2}a\right)\Gamma\left(\frac{1}{2}a-c+\frac{1}{2}\right)\Gamma\left(\frac{1}{2}a-b+\frac{3}{2}\right)\Gamma(2+a-b-c)}$$

$$+ \frac{\beta}{(b-1)} \frac{2^{-a-1}\,\Gamma\left(\frac{1}{2}\right)\Gamma(1+a-c)\,\Gamma(1+a-b)\,\Gamma\left(1+\frac{1}{2}a-b-c\right)}{\Gamma\left(\frac{1}{2}a+\frac{1}{2}\right)\Gamma\left(\frac{1}{2}a-b+1\right)\Gamma\left(\frac{1}{2}a-c+1\right)\Gamma(1+a-b-c)} \qquad (1.9)$$

provided $\Re(a-2b-2c) > -2$ and $\alpha$ and $\beta$ are given by

$$\alpha = 1 - \frac{1}{d}(1+a-b) \qquad (1.10)$$

$$\beta = \frac{1+a-b}{1+a-b-c}\left[\frac{a}{d}(1+a-b-2c) - 2\left(\frac{a}{2}-b-c+1\right)\right] \qquad (1.11)$$

Clearly for $d = c$, (1.8) reduces to (1.6) and for $d = 1 + a - b$, (1.9) reduces to (1.7) respectively.

Applications of the above mentioned classical summation theorems are well known now. It has been pointed out by Berndt [1] that several interesting summations due to Ramanujan can be obtained quite simply by employing above mentioned classical summation theorems.

Here, we would like to mention the following interesting summations due to Ramanujan[1,5,6]. These are

$$1 + \frac{1}{5}\left(\frac{1}{2}\right)^2 + \frac{1}{9}\left(\frac{1.3}{2.4}\right)^2 + \cdots = \frac{\pi^2}{4\,\Gamma^4\left(\frac{3}{4}\right)} \tag{1.12}$$

$$1 + \frac{1}{5^2}\left(\frac{1}{2}\right) + \frac{1}{9^2}\left(\frac{1.3}{2.4}\right) + \cdots = \frac{\pi^{\frac{5}{2}}}{8\sqrt{2}\,\Gamma^2\left(\frac{3}{4}\right)} \tag{1.13}$$

and

$$1 + \frac{1}{5}\left(\frac{1}{2}\right) + \frac{1}{9}\left(\frac{1.3}{2.4}\right) + \cdots = \frac{\pi^{\frac{3}{2}}}{2\sqrt{2}\,\Gamma^2\left(\frac{3}{4}\right)} \tag{1.14}$$

As explained by Berndt[1] that the summations (1.12) and (1.13) due to Ramanujan can be obtained quite simply by employing classical Dixon summation theorem (1.7) by taking

(i) $\quad a = b = \frac{1}{2}$, $c = \frac{1}{4}$ and

(ii) $\quad a = \frac{1}{2}$, $b = c = \frac{1}{4}$

respectively.

The summation (1.14) obtained by Ramanujan[6] using an integral representation.

Very recently Rathie and Paris [8] pointed out that the summation (1.14) can also be obtained in a very simple manner by employing Gauss summation theorem (1.6) by taking $a = \frac{1}{2}$, $b = \frac{1}{4}$ and $c = \frac{5}{4}$.

In 2010, Kim, et al. [2] obtained the extensions of Ramanujan summations (1.12) and (1.13) by employing extension of Dixon summation theorem (1.9) in the following form {(1.15) written here with minor correction}

$$1 + \frac{1}{5}\left(\frac{1}{2}\right)^2 \left(\frac{d+1}{2d}\right) + \frac{1}{9}\left(\frac{1.3}{2.4}\right)^2 \left(\frac{d+1}{3d}\right) + \cdots = \frac{4}{3\pi}\left(\frac{1}{d} - 1\right) + \frac{\pi^2}{3\Gamma^4\left(\frac{3}{4}\right)}\left(1 - \frac{1}{4d}\right) \quad (1.15)$$

and

$$1 + \frac{1}{5^2}\left(\frac{1}{2}\right)\left(\frac{5(d+1)}{9d}\right) + \frac{1}{9^2}\left(\frac{1.3}{2.4}\right)\left(\frac{5(d+2)}{13d}\right) + \cdots$$

$$= \frac{5\pi^{\frac{3}{2}}}{48\sqrt{2}\,\Gamma^2\left(\frac{3}{4}\right)}\left(\frac{5}{4d} - 1\right) - \frac{5}{32\sqrt{2}}\frac{\pi^{\frac{5}{2}}}{\Gamma^2\left(\frac{3}{4}\right)}\left(\frac{1}{4d} - 1\right) \quad (1.16)$$

each for $d \neq 0, -1, -2, \cdots$

In this short research note, we aim to give an extension of Ramanujan summation (1.14).

## 2. Extension of Ramanujan summation (1.14):

The summation to be established in this note is

$$1 + \left(\frac{1}{2}\right)\left(\frac{1}{9}\right)\left(\frac{d+1}{d}\right) + \left(\frac{1.3}{2.4}\right)\left(\frac{1.5}{9.13}\right)\left(\frac{d+2}{d}\right) + \cdots = \frac{5\pi^{\frac{3}{2}}}{12\sqrt{2}\,\Gamma^2\left(\frac{3}{4}\right)}\left(1 + \frac{1}{4d}\right) \quad (2.1)$$

for $d \neq 0, -1, -2, \cdots$

**Proof :** In order to derive (2.1), we proceed as follows. Setting $a = \frac{1}{2}$, $b = \frac{1}{4}$ and $c = \frac{5}{4}$ in extended Gauss summation theorem (1.8), we see that

$$L.H.S. = {}_3F_2\left[\begin{array}{c}\frac{1}{2},\ \frac{1}{4},\ d+1 \\ \frac{9}{4},\ d\end{array};1\right]$$

Expressing ${}_3F_2$ as a series with the help of the definition (1.1), we have

$$= \sum_{n=0}^{\infty} \frac{\left(\frac{1}{2}\right)_n \left(\frac{1}{4}\right)_n (d+1)_n}{\left(\frac{9}{4}\right)_n (d)_n\, n!}$$

Using Pochhammer symbol and expanding it, we have after some simplification

$$= 1 + \left(\frac{1}{2}\right)\left(\frac{1}{9}\right)\left(\frac{d+1}{d}\right) + \left(\frac{1.3}{2.4}\right)\left(\frac{1.5}{9.13}\right)\left(\frac{d+2}{d}\right) + \cdots \tag{2.2}$$

Also with the same substitution in the right-hand side, we have

$$R.H.S. = \frac{1}{2}\frac{\Gamma\left(\frac{9}{4}\right)\Gamma\left(\frac{1}{2}\right)}{\Gamma\left(\frac{7}{4}\right)}\left(1 + \frac{1}{4d}\right)$$

Using the results $\Gamma(z+1) = z\,\Gamma(z)$, $\Gamma(z)\Gamma(1-z) = \frac{\pi}{sin\pi z}$, and $\Gamma\left(\frac{1}{2}\right) = \sqrt{\pi}$, we

have after some simplification,

$$= \frac{5\pi^{\frac{3}{2}}}{12\sqrt{2}\,\Gamma^2\left(\frac{3}{4}\right)}\left(1 + \frac{1}{4d}\right) \tag{2.3}$$

Finally equating (2.2) and (2.3), we get (2.1). This completes the proof of (2.1).

We conclude this note with the remark that our extended summation (2.1) reduces to Ramanujan summation by taking $d = \frac{5}{4}$. Moreover, by giving different values for d in our general summation (2.1), we can obtain a large number of new summations.